# Integer roots of quadratic and cubic polynomials with integer coefficients

Konstantine Zelator

Mathematics, Computer Science and Statistics

212 Ben Franklin Hall

Bloomsburg University

400 East Second Street

Bloomsburg, PA  17815

U.S.A.

kzelator@bloomu.edu

konstantine_zelator@yahoo.com

Konstantine Zelator

P.O. Box 4280

Pittsburgh, PA  15203

# 1. INTRODUCTION

The subject matter of this work is quadratic and cubic polynomials with integral coefficients; which also has all of the roots being integers. The purpose of this work is to determine precise (i.e. necessary and sufficient) coefficient conditions in order that a quadratic or a cubic polynomial have integer roots only.

The results of this paper are expressed in Theorems 3, 4, and 5.

The level of the material in this article is such that a good second or third year mathematics major; with some exposure to number theory (especially the early part of an introductory course in elementary number theory); can comfortably come to terms with.

Let us outline the organization of this article. There are three lemmas and five theorems in total. The results expressed in Theorems 4 and 5; are not found in standard undergraduate texts (in the United States) covering material of the first two years of the undergraduate mathematics curriculum. But some of these results might be found in more obscure analogous books (probably out of print) around the globe.

The three lemmas are number theory lemmas and can be found in Section 3. Lemma 1 is known as Euclid's lemma; it is an extremely well known lemma and can be found in pretty much every introductory book in elementary number theory. We use Lemma 1 to establish Lemma 3; which is in turn used in the first proof of Theorem 3 (we offer two proofs for Th. 3). Lemma 2, also known as the nth power lemma, is also very well known; but perhaps a bit less than Lemma 1. A reference for Lemma 1 can be found in [1].

A reference for Lemma 2 can be found in [2] (it is also listed as an exercise in [1]). Lemma 2 is used in the second proof of Theorem 3 and also toward establishing Theorem 4. Theorem 3, in Section 4; states necessary and sufficient conditions for a quadratic polynomial with integer coefficients to have two integral roots or zeros. The result in Th. 3 is generally recognizable, but the interested reader may not find it in standard undergraduate texts. We offer two proofs to Theorem 3.

Theorem 4, in Section 6; is the one that (most likely) the reader will be the least familiar with. Theorem 4 gives precise (i.e. necessary and sufficient) coefficient conditions, in order that a monic (i.e. with leading coefficient 1) cubic polynomial with integer coefficients; have one double (i.e. multiplicity 2) integral root, as well as another integer root.

In Section 7, Theorem 5, parametrically describes a certain family of quadratic polynomials $ax^2 + bx + c$ with integer coefficients; and with integral zeros. As for Theorems 1 and 2, they are extremely well known; and both found in Section 2. These two theorems; at least in part, can be found pretty much in every college algebra/pre-calculus book (for example see reference [3]); not to mention more advanced books. Theorem 1 combines two very well-known results into one. It states that every nth degree polynomial with complex number coefficients; has exactly n complex roots, counting multiplicities. And it (Th. 1) also gives the factorization of a polynomial into $n$ linear factors. See reference [4].

Lastly, Theorem 2 is the familiar rational root theorem.
We start with two definitions.
*Definition 1: A proper rational number is a rational number which is not an integer.*
*Definition 2. A monic polynomial (with complex coefficients), is a polynomial whose leading coefficient is equal to 1.*

## 2. SOME PRELIMINARIES FROM ALGEBRA

First, we state the very well-known Theorems 1 and 2.

**Theorem 1** (n roots and linear factorization theorem)
*Let $p(x) = a_n x^n + a_{n-1} x^{n-1} + \cdots + a_1 x + a_0$; be a polynomial function of degree $n \geq 1$, and so $a_n \neq 0$; where the coefficients $a_n, a_{n-1}, \ldots, a_1 a_0$ are complex numbers.*
*(a) The polynomial $p(x)$ has exactly $n$ complex roots, counting multiplicities*
*(b) $p(x) = a_n(x - r_1)(x - r_2) \ldots (x - r_n)$, where $r_1, r_2, \ldots, r_n$ are the roots of $p(x)$.*

In more advanced algebra courses; students learn that a polynomial function (in one variable $x$) over an algebraically closed field $F$; and with coefficients from the field $F$; has exactly $n$ zeros in $F$, counting multiplicities. Such is the case, in particular, when $F = \mathbb{C}$, the field of complex numbers.

For more details, see reference [4].

For Theorem 2 below, see reference [1] and [3]

**Theorem 2** (Rational Root Theorem)
*Let $p(x) = a_n x^n + a_{n-1} x^{n-1} + \cdots + a_1 x + a_0$; of degree $n \geq 1$, and so $a_n \neq 0$; with integer coefficients $a_n, a_{n-1}, \ldots, a_1, a_0$.*

*Suppose that the rational number $r = \dfrac{k}{\ell}$, in lowest terms (i.e. $k$ and $\ell$ are*

*relatively prime integers, with $\ell \neq 0$); is a zero or root of $p(x)$. The, the integer $k$*
*is a divisor of the consant term $a_o$; adn the integer $\ell$ is a divisor of the leading coefficient $a_n$.*

---

It follows immediately from Th. 1(b) (with a straightforward expansion/calculation), that if $t(x) = a_2 x^2 + a_1 x + a_0$ is a quadratic trinomial with complex number coefficients $a_2, a_1, a_0$; and roots $r_1$ and $r_2$. Then $t(x) = a_2(x - r_1)(x - r_2)$; and consequently $r_1 + r_2 = -\frac{a_1}{a_2}$ and $r_1 r_2 = \frac{a_0}{a_2}$

(1)

---

Likewise, if $c(x) = a_3 x^3 + a_2 x^2 + a_1 x + a_0$ is a cubic polynomial function with complex number coefficients $a_3, a_2, a_1, a_0$; and roots $r_1, r_2, r_3$. Then Theorem 1, part (b) implies that

$$c(x) = a_3(x - r_1)(x - r_2)(x - r_3); \text{ and accordingly,}$$

$$\begin{cases} r_1 + r_2 + r_3 = -\frac{a_2}{a_3} \\ r_1 r_2 + r_2 r_3 + r_3 r_1 = \frac{a_1}{a_3} \\ r_1 r_2 r_3 = -\frac{a_1}{a_3} \end{cases}$$

(2)

---

In the case of the trinomial $t(x) = a_2 x^2 + a_1 x + a_0$; we also have the familiar quadratic formula that gives us the two roots $r_1$ and $r_2$:

$$r_1 = \frac{-a_1 + \sqrt{a_1^2 - 4a_2 a_0}}{2a_2}, \quad r_2 = \frac{-a_1 - \sqrt{a_1^2 - 4a_2 a_0}}{2a_2}$$

(3)

### 3. THREE LEMMAS FROM NUMBER THEORY

**Lemma 1** (Euclid's lemma) (see [1])
*Suppose that a, b, c, are positive integers, and that a is a divisor of the product bc. Also, assume that the integer a is relatively prime to the integer c. Then, a is a divisor of c.*

**Lemma 2** (nth power lemma) (see [2])
*Let n be a positive integer. Then, a positive integer b is the nth power of a positive rational number; if, and only if, b is the nth power of a positive integer. In particular, b is a rational square; if, and only if it is an integer square.*

**Lemma 3**
Suppose that the sum and the product of two rational numbers; are both integers. Then, both of these two rational numbers are integers.

**Proof:** Let $r_1$, $r_2$ be two such rational numbers.

Then, $\left\{\begin{array}{l} r_1 + r_2 = i_1 \\ r_1 r_2 = i_2 \end{array}\right\}$ (4)

Where $i_1, i_2$ are integers.

If $r_1 = 0$, then $r_2 = i_1$, and we are done. Similarly if $r_2 = 0$. Next assume that both $r_1$ and $r_2$ are nonzero. We write $r_1$ and $r_2$ in lowest terms with positive denominators:

$$\left\{\begin{array}{c} r_1 = \frac{k_1}{\ell_1} , r_2 = \frac{k_2}{\ell_2}; \text{ where } k_1, k_2 \text{ are integers;} \\ \text{and } \ell_1, \ell_2 \text{ are positive integers such that} \\ (k_1, \ell_1) = 1 = (k_2, \ell_2); \text{ and so} \\ (|k_1|, \ell_1) = 1 = (|k_2|, \ell_2) \text{ as well} \end{array}\right\}$$ (5a)

Note: we have used the standard notation $(\alpha, \beta)$; to denote the greatest common divisor between two integers $\alpha$ and $\beta$.

Combining (4) with (5a), we obtain

$$\left\{\begin{array}{c} k_1 \ell_2 + k_2 \ell_1 = i_1 \ell_1 \ell_2 \\ k_1 k_2 = i_2 \ell_1 \ell_2 \end{array}\right\}$$ (5b)

The first equation in (5b) implies,
$k_1 \ell_2 = \ell_1 \cdot (i_1 \ell_2 - k_2)$; which in turn implies
$$|k_1| \ell_2 = \ell_1 |i_1 \ell_2 - k_2|$$ (5c)

According to (5c), the positive integer $\ell_1$ is divisor of the product $|k_1| \ell_2$; and since by (5a), it is relatively prime to $|k_1|$; it follows by Lemma 1, that $\ell_1$ is a divisor of the positive integer $\ell_2$. A similar argument (using once again the first equation in (5b)); shows that $\ell_2$ is a divisor $\ell_1$ as well. Since $\ell_1$ and $\ell_2$ are positive integers which are divisors of each other; they must be equal: $\ell_1 = \ell_2$. Hence by the second equation in (5b) we get,
$$k_1 k_2 = i_2 \ell_1^2$$ (5d)

Since $\ell_1 = \ell_2$; (5a) implies $((k_1, \ell_1) = 1 = (k_2, \ell_1)$; which in turn gives $(k_1 k_2, \ell_1) = 1$ (an exercise in elementary number theory). This last condition, in conjunction with (5d), implies that $\ell_1^2 = 1$; and so $\ell_1 = 1 = \ell_2$. Hence $r_1 = k_1, r_2 = k_2$; both $r_1$ and $r_2$ are integers.

# 4.  THEOREM 3 AND ITS TWO PROOFS

**Theorem 3**

*Let $q(x) = ax^2 + bx + c$, being a quadratic trinomial with integer coefficients a, b, c. Then, both roots or zeros of $q(x)$ are integers; if, and only if,*

*(i)     The integer $b^2 - 4ac$ is an integer or perfect square. And (ii) The leading coefficient a is a divisor of both b and c.*

***First proof:*** First assume that both roots $r_1$ and $r_2$ are integers. From (1) we have,

$r_1 + r_2 = -\frac{b}{a}$ and $r_1 r_2 = \frac{c}{a}$ ; Which clearly implies that both $\frac{b}{a}$ and $\frac{c}{a}$ are integers; which in turn implies that a is a divisor of both coefficients b and c.

By completing the square we have,

$$q(x) = a(x + \frac{b}{2a})^2 + \frac{4ac - b^2}{4a}$$

Since $r_1$ is a root; we have $q(r_1) = 0$. And so,

$a(r_1 + \frac{b}{2a})^2 + \frac{4ac-b^2}{4a} = 0$; or equivalently (multiply both sides by 4a)

$(2ar_1 + b)^2 = b^2 - 4ac$; which proves that $b^2 - 4ac$ is a perfect or integer square, since $r_1$ is an integer (and thus, so is $2ar_1 + b$).

Now the converse.

Assume that $b^2 - 4ac = k^2$; where k is a nonnegative integer. And that a is a divisor of both b and c. From the quadratic formula (3), we have

$r_1 = \frac{-b+k}{2a}$ and $r_1 = \frac{-b-k}{2a}$; which shows that both $r_1$ and $r_2$ are both rational numbers. Moreover, $\frac{b}{a}$ and $\frac{c}{a}$ are both integers (since a is a divisor of b and c). Thus, both the sum $r_1 + r_2 = \frac{b}{a}$; and the product $r_1 r_2 = \frac{c}{a}$ are integers. Hence by Lemma 3 it follows that both $r_1$ and $r_2$ are integers.

**Second Proof:** First, assume that both $r_1$ and $r_2$ are integers. As in the first proof; it follows from $r_1 + r_2 = -\frac{b}{a}$ and $r_1 r_2 = \frac{c}{a}$; that a is a divisor of both b and c. Furthermore from the quadratic formula (3); we have that $r_1 = \frac{-b+\sqrt{b^2-4ac}}{2a}$. Which implies that since $r_1$, b, a are all rationals; the number $\sqrt{b^2 - 4ac}$; must also be rational.

That is, the integer $b^2 - 4ac$ must be a rational square. Hence; by Lemma 2; the integer $b^2 - 4ac$ must be a perfect or integer square.

Let us prove the converse. Suppose that a is a divisor of both b and c. And also, assume that $b^2 - 4ac$ is a perfect square. We have,

$$\left\{ \begin{array}{c} b = a \cdot u, c = a \cdot v, b^2 - 4ac = m^2; \\ where\ u, v, m\ are\ integers\ with\ m \geq 0 \end{array} \right\} \tag{6}$$

From (6) we get, $a^2(u^2 - 4v) = m^2$ \hfill (7)

Equation (7) shows that $a^2$ is a divisor of $m^2$; and therefore a must be a divisor of m (an exercise in elementary number theory). We put,

$$m = a \cdot t, t\ an\ integer \tag{8}$$

Combining (8) with the quadratic formula (3); we see that,

$$r_1 = \frac{a \cdot u + |a||t|}{2a}\ and\ r_2 = \frac{-a \cdot u - |a||t|}{2a}$$

If $a > 0$ we obtain,

$$r_1 = \frac{u + |t|}{2}\ and\ r_2 = \frac{-(u + |t|)}{2} \tag{9i}$$

While if $a < 0$, we have,

$$r_1 = \frac{u - |t|}{2}\ and\ r_2 = \frac{-u + |t|}{2} \tag{9ii}$$

Also, combining (7) with (8) yields,

$$u^2 - 4v = t^2 \tag{10}$$

Equation (10) clearly shows that the integers u and t have the same parity: they are both even, or both odd. Hence (9i) and (9ii) imply that $r_1\ and\ r_2$ are both integers.

## 5. FOUR OBSERVATIONS

### Observation 1

*If an nth degree polynomial p(x) has integer leading coefficient $a_n$; and if all the roots $r_1, r_2, ..., r_n$; are also integers. Then, it follows from Theorem 1, part (b); that all the other coefficients $a_{n-1}, ..., a_1, a_0$ must also be integers (just expand the right-hand side and collect terms); and that the leading coefficient $a_n$ must be a divisor of each one of the other coefficients* $a_{n-1}, \ldots, a_1, a_0$. *And so the corresponding monic polynomial,* $\frac{1}{a_n} \cdot p(x) = x^n + \frac{a_{n-1}}{a_n} \cdot x^{n-1} + \cdots + \frac{a_1}{a_n} + \frac{a_0}{a_n}$; *has integer coefficients, and of course, the exact same integer roots as p(x).*

**Observation 2**

*If p(x) is an nth degree monic polynomial with integral coefficients, that is $a_n = 1$ and $a_{n-1}, \ldots, a_1, a_0 \in \mathbb{Z}$. Then, a consequence of Theorem 2 is that if p(x) has at least one rational root. Then, each of its rational roots or zeros, much be an integer. Equivalently, p(x) cannot have a rational root which is a proper rational number. Further equivalently, each real root (if any) of p(x); must be either an integer or otherwise an irrational number.*

**Observation 3**

*Note that Observation 2 does imply Lemma 3. Inded, suppose that $r_1$ and $r_2$ are rational numbers such that $r_1 + r_2 = \alpha$ and $r_1 r_2 = \beta$; were $\alpha, \beta$ are inegers. Obviously $r_1$ and $r_2$ are the two rational roots of the monic quadratic trinomial $t(x) = (x - r_1)(x - r_2)$;*

$t(x) = x^2 - (r_1 + r_2)x + r_1 r_2$;

$t(x) = x^2 - \alpha x + \beta$

*Since t(x) is a monic polynomial with integer coefficients; each of its rational roots must be integral, by Observation 2. Hence, both $r_1$ and $r_2$ must be integers.*

**Observation 4**

An immediate consequence of Theorem 3 is that the monic trinomial $x^2 + bx + c$ (with b, c integers) has two integers roots; if and only if $b^2 - 4c$ is an integral square. Similarly, the trinomial $-x^2 + bx + c$ has two integer roots; if and only if the integer $b^2 + 4c$ is the square of an integer.

## 6. FIVE SPECIAL CASES OF CUBICS

In this section, we consider five special cases of cubic or degree 3 monic polynomials with integer coefficients;

$$c(x) = x^3 + bx^2 + cx + d; b, c, d \in \mathbb{Z}$$

And with all three roots of c(x) being integers

*Special Case 1: The number 0 is a root of c(x); so that d=0*

Let $r_1$ and $r_2$ be the other two integer roots; and $r_3 = 0$. By (2) we have,

$$r_1 + r_2 = -b$$
$$r_1 r_2 = c$$

Clearly $r_1$ and $r_2$ are the two integer roots of the monic trinomial

$$(x - r_1)(x - r_2) = x^2 - (r_1 + r_2)x + r_1 r_2; \ x^2 + bx + c$$

By Observation 4, it follows that this will happen precisely when $b^2 - 4ac$ is an integer square.

*Summary*

Let $c(x) = x^3 + bx^2 + cx + d$ be a monic cubic polynomial with integer coefficients. Then $c(x)$ has three integer roots with zero being one of them; precisely when the following conditions are satisfied;
$d = 0$ and $b^2 - 4c = k^2$; k a nonnegative integer
Then, the three integer roots of $c(x)$ are:
$$r_1 = \frac{-b + k}{2}, r_2 = \frac{-(b + k)}{2}, r_3 = 0$$

*Special Case 2: The number 1 is a root of c(x)*

Clearly, this occurs if, and only if, the condition $b + c + d = 1$ is satisfied.

Let $r_1, r_2$, and $r_3 = 1$, be the three integer roots of *c(x)*. From (2) we get

$$\begin{Bmatrix} r_1 + r_2 = -b - 1 \\ r_1 + r_2 + r_1 r_2 = c \\ r_1 r_2 = -d \end{Bmatrix}$$

Note that the second equation results by adding the first and third equations memberwise; due to the condition $b + c + d = 1$. The numbers $r_1$ and $r_2$ are the two roots of the quadratic trinomial $(x - r_1)(x - r_2) = x^2 + (b - 1)x - d$. Once again, Observation 4 tells us that is quadratic trinomial will have two integer roots if, and only if, $(b + 1)^2 + 4d$ is an integer square

*Summary*

*Let $c(x) = x^3 + bx^2 + cx + d$ be a monic cubic polynomial with integer coefficients. Then c(x) has three integer roots with 1 being among them; precisely if the following conditions are satisfied:*

$b + c + d = -1$ $and$ $(b + 1)^2 + 4d = k^2$, $k$ $a$ $nonnegative$ $integer.$ $Then,$ $the$ $three$ $roots$
$are:$ $r_1 = \frac{-(b+1)+k}{2}, r_2 = \frac{-(b+1+k)}{2}, r_3 = 1$

*Special Case 3: The number -1 is a root of c(x)*

This case is analogous to Special Case 2. We invite the reader to fill in the details.

*Summary*

Let $c(x) = x^3 + bx^2 + cx + d$ be a monic cubic polynomial with integer coefficients. Then c(x) has three integer roots with -1 being among them; precisely when $b - c + d = 1$ $and$
$(b - 1)^2 - 4d = k^2$; k a nonnegative integer. Then, the three roots are:
$$r_1 = \frac{b-1+k}{2}, r_2 = \frac{b-1-k}{2}, r_3 = -1$$

*Special Case 4: c(x) has a triple (i.e. multiplicity 3) integer root*

Let r be that root. We have,

$$c(x) = (x - r)^3 = x^3 + bx + c + d;$$
$$x^3 - 3rx^2 + 3r^2x - r^3 = x^3 + bx + c + d$$

(Identically equal as polynomial functions);

Which $\qquad -3r = b, 3r^2 = c, and -r^3 = d;$
$$r = -\frac{b}{3}, c = \frac{b^2}{3}, d = \frac{b^3}{27}$$

*Summary*

Let $c(x) = x^3 + bx^2 + cx + d$ be a monic polynomial with integer coefficients. Then c(x) has a triple integer root r if, and only if,
$$c = \frac{b^2}{3} \text{ and } d = \frac{b^3}{27}$$
The integer root is $r = -\frac{b}{3}$

*Special Case 5: c(x) has a double integer root and another integer root*

Let M be the integer root of multiplicity 2; and N, distinct from M; the integer root of multiplicity 1. First, by (2) we obtain

$$\begin{cases} 2M + N = -b \\ M^2 + 2MN = c \\ M^2N = -d \end{cases} \qquad \begin{matrix} (10i) \\ (10ii) \\ (10iii) \end{matrix}$$

Now, $(x) = x^3 + bx^2 + cx + d$ ; as we know from calculus, its derivative function is $c'(x) = 3x^2 + 2bx + c$. Also, from theorem 1, part (b); we have $c(x) = (x - M)^2 \cdot (x - N)$. Putting these three together we have,

$$\begin{cases} c(x) = x^3 + bx^2 + cx + d = (x - N)(x - M)^2 \\ \qquad and \ c'(x) = 3x^2 + 2bx + c \end{cases} \tag{11}$$

Also from the product rule in calculus we have,

$$c'(x) = (x - M)^2 + 2(x - N)(x - M); and \ so$$

$$c'(M) = 0; the \ integer \ M \ is \ a \ zero \ or \ root \ of \ the \ quadratic \ trinomial \ c'(x) as \ well$$

The coefficients of c'(x) are the integers 3,2b, and c. Since M is an integer; it then follows from (1), that the other root of the trinomial c'(x); must be a rational number. (See note below)

Hence by (3) it follows that the real number,

$\sqrt{(2b)^2 - 4(3)c} = \sqrt{4(b^2 - 3c)}$; must be a rational number; which in turn implies that the integer $4(b^2 - 3c)$ is the square of rational number. Hence, by Lemma 2; the integer $4(b^2 - 3c)$ must be an integer square. We have, $4(b^2 - 3c) = K^2$; K a nonnegative integer.

Clearly k must be even; K= $2k$; for some nonnegative integer k. Thus,

$$k \geq 0, K = 2k \ and \ b^2 - 3c = k^2 \tag{12}$$

Note: Above, one can argue without using Lemma 2. Simply by saying that since M is a root. Then by (3),
$$\sqrt{4(b^2 - 3c)} = 6M + 2b; or \ alternatively \sqrt{4(b^2 - 3c)} = -(6M + 2b).$$
In either case, $4(b^2 - 3c) = (6M + 2b)^2; an \ integer \ square$

Let us continue. Equation (12) clearly shows that either both integers b and k are divisible by 3; or neither of them is. So, we consider two cases. The case in which neither b or k is a multiple of 3; and the case with both b and k being multiples of 3.

**Case 1:** $b \not\equiv 0 (mod \ 3) and \ so \ k \ \not\equiv 0 \ (mod \ 3) as \ well.$

The two roots of the quadratic trinomial $c'(x) = 3x^2 + 2bx + c \ are \ \frac{-2b+k}{6} \ and \ \frac{-(2b+k)}{6}$. One of these two rational numbers is the integer M. Note that k must be even; otherwise (if k were odd) both of these numbers would be proper rationals. Furthermore, if $b \equiv k (mod \ 3)$ (i.e. both b and k congruent to 1(mod3); or both congruent to 2(mod3)); then $M = \frac{-(2b+k)}{6}$, while the other root $\frac{-2b+k}{6}$ is a proper rational. On the other hand, if $b \not\equiv k (mod \ 3)$ (i.e. one of b, k; is congruent to 1(mod3), the other congruent to 2(mod3), then $M = \frac{-2b+k}{6}$; while $\frac{-(2b+k)}{6}$ is a proper rational.

The value of N is determined by $(10i)$; and the value of d is determined by $(10iii)$. A routine calculation shows that the equation $(10ii)$ is also satisfied (by using (12) as well). Note that k actually cannot equal zero; for this would imply $M = -\frac{b}{3}$, a triple root of the cubic polynomial $c(x)$; contrary to the assumption that M is a double root.

**Case 2:**  $b \equiv k \equiv 0(mod3)$

Since both b and k are divisible by 3; equation (12) implies that c must also be divisible by 3 as well. We put,

$$\left\{\begin{array}{c} b = 3B, c = 3C, k = 3t; \\ B, C, t \in \mathbb{Z} \end{array}\right\} (13)$$

And thus by (12) and (13) we get,

$$B^2 - C = t^2 \ (14)$$

The quadratic trinomial $c'(x)$ takes the form,

$$c'(x) = 3(x^2 + 2Bx + C)$$

Clearly then, M is one of the two integer roots of the monic quadratic trinomial $x^2 + 2Bx + C$. These two integers are $-B + t$ and $-(B + t)$. Again, as in Case 1; the value of N is determined by $(10i)$; and the value of d by $(10iii)$. And, one can easily check that $(10ii)$ is also met. And t cannot be zero, since M is a double root; not a triple root.

We summarize the results of Cases 1 and 2 in Theorem 4 below.

**Theorem 4**

*Consider the monic cubic polynomial with integer coefficients, $c(x) = x^3 + bx^2 + cx + d$. Then, $c(x)$ has an integer root M of multiplicity 2, and an integer root N of multiplicity 1. If and only if, the coefficients b, c, d; satisfy precisely one of the group conditions:*

*Group 1:  $b \not\equiv 0(mod3); k \not\equiv 0(mod3), k$ an even positive integer, $b \equiv k(mod3), c = \frac{b^2-k^2}{3}$ (note that $b^2 \equiv k^2 \equiv 1(mod3)$; since $bk \not\equiv 0(mod3)$), and $d = (\frac{2b+k}{6})^2 \cdot \left[b - \left(\frac{2b+k}{3}\right)\right]$. The roots are $M = -\left(\frac{2b+k}{6}\right)$ and $N = -b + \frac{(2b+k)}{3}$.*

*Group 2:  $b \not\equiv 0(mod3); k \not\equiv 0(mod3), k$ an even positive integer, $b \not\equiv k(mod3), c = \frac{b^2-k^2}{3}$, and $d = (\frac{-2b+k}{6})^2 \cdot \left[b + \frac{-2b+k}{3}\right]$. The roots are $M = \frac{-2b+k}{6}$ and  $N = -b + \frac{2b-k}{3}$.*

*Group 3:  $b = 3B, c = 3C, C = B^2 - t^2$; with B, C being integers; and t a positive integer.  Also with $d = (-B + t)^2 \cdot [B + 2t]$.  The roots are $M = -(B + t)$ and $N = -(B + 2t)$.*

*Group 4:  $b = 3B, c = 3C, C = B^2 - t^2$; with B, C being a integers; and t a positive integer. Also with $d = (B + t)^2 \cdot [B - 2t]$.*

*The roots are $M = -(B + t)$ and $N = -B + 2t$.*

## 7.  A FAMILY OF QUADRATIC TRINOMIALS WITH INTEGER COEFFICIENTS AND INTEGRAL ROOTS

In this last section we endeavor to describe a certain type of quadratic trinomials $t(x) = ax^2 + bx + c$.

*These are trinomials with two integer roots and with the integer coefficients a, b, c; satisfying $|a| > 1$ and $|b| = p$; where p is a prime number.*
*Analysis of the above type of quadratic trinomials*

By Theorem 3, since $t(x)$ has two integer roots $r_1$ and $r_2$; we must have (precise conditions),

$$\left\{ \begin{array}{c} b^2 - 4ac = k^2 \text{ and both b and c} \\ \text{being divisible by the integer } a; k \text{ a nonnegative integer} \end{array} \right\} \tag{15i}$$

Since $|b| = p$, a prime; and since b is divisible by a; it easily follows that  $|a|$=1 or p. But $|a| > 1$; thus we must have  $|a| = p$; which means that $a = p \ or - p$.  Also since c is a multiple of a; we have

$$v \cdot a = c \quad (15ii)$$

We have $a^2 = p^2$, then $(15i)$ and $(15ii)$ yield, since  $|b|$ = p;

$$p^2(1 - 4v) = k^2,$$

Which implies that p must be a divisor of k; $k = pT, T$ a positive integer. And we obtain, from the last equation,

$$v = \frac{1 - T^2}{4}; T \text{ an odd positive integer } \left( so \text{ that } T^2 \equiv 1 (mod 4) \right).$$

*Definition 3:  A quadratic trinomial $t(x) = ax^2 + bx + c$ with integer coefficients is said to be a p-Type 1; where p is a prime; if it has two integer roots and the coefficients satisfy the conditions $|a| > 1$ and $|b| = p$.*

We can now state Theorem 5.

**Theorem 5**

*For a given prime p, there are exactly four groups p-Type 1 quadratic trinomials with integer coefficients.*

*Group 1:  Coefficients are $a = p, b = p, c = \frac{P(1-T^2)}{4}$; where T is an odd positive integer.  The two integer roots are $r_1 = \frac{-1+T}{2}$, $r_2 = \frac{-(1+T)}{2}$.*

*Group 2:  Coefficients are $a = -p, b = p, c = \frac{P(1-T^2)}{4}$; where T is an odd positive integer.  The two integer roots are $r_1 = \frac{1-T}{2}$, $r_2 = \frac{1+T}{2}$.*

*Group 3:  Coefficients are $a = -p, b = -p, c = \frac{P(1-T^2)}{4}$; where T is an odd positive integer. The two integer roots are $r_1 = \frac{-(1+T)}{2}$, $r_2 = \frac{-1+T}{2}$.*

*Group 4:  Coefficients are $a = p, b = -p, c = \frac{P(1-T^2)}{4}$; where T is an odd positive integer.  The two integer roots are $r_1 = \frac{1+T}{2}$, $r_2 = \frac{1-T}{2}$.*

Note: We see that each of the above four groups; is a one-parameter family. The parameter being the odd positive integer T.